\documentclass{article}

\usepackage{latexsym}
\input{amssym}

\newcommand{\iso}{\cong}
\newcommand{\sub}{\subseteq}
\newcommand{\nsub}{\trianglelefteq}
\newcommand{\zed}{\Bbb{Z}}
\newcommand{\pf}{\paragraph*{Proof}}
\newcommand{\done}{\hfill\ensuremath{\Box}}
\newcommand{\defemph}[1]{\textbf{#1}}
\newcommand{\sect}[1]{\section{#1}\setcounter{paragraph}{0}}
\newcommand{\go}{\longrightarrow}
\newcommand{\of}{\raisebox{0.2mm}{\ensuremath{\scriptstyle\circ}}}

\title{Perfect Numbers and Groups}
\author{Tom Leinster\\ \\
	\normalsize{Department of Pure Mathematics, University of
	Cambridge}\\ 
	\normalsize{Email: leinster@dpmms.cam.ac.uk}\\
	\normalsize{Web: http://www.dpmms.cam.ac.uk/$\sim$leinster}}
\date{1st April, 2001}

\begin{document}


\begin{titlepage} 

\maketitle
\thispagestyle{empty}

\vfill

\begin{center}	\bfseries
Abstract
\end{center}
\begin{quotation}
A number is perfect if it is the sum of its proper divisors; here we call a
finite group `perfect' if its order is the sum of the orders of its proper
normal subgroups.  (This conflicts with standard terminology but confusion
should not arise.)  The notion of perfect group generalizes that of perfect
number, since a cyclic group is perfect just when its order is perfect.  We
show that, in fact, the only abelian perfect groups are the cyclic ones, and
exhibit some non-abelian perfect groups of even order.
\end{quotation}

\vfill

\begin{quotation}	\emph{%
This article was originally composed in~1996 for \emph{Eureka}, the journal
of the Cambridge student mathematical society (but has yet to appear, as
no issue has been published since).  It is therefore written to be
comprehensible to an undergraduate readership, and contains many reminders of
basic facts.}
\end{quotation}

\vfill

\begin{center}
\bfseries Contents
\end{center} 

\begin{center}
\begin{tabular}{rlcr}
\ref{sec:numbers}	&Perfect numbers	&&\pageref{sec:numbers}	\\
\ref{sec:basics} 	&Definition and first examples of perfect groups
						&&\pageref{sec:basics} 	\\
\ref{sec:mult}		&Multiplicativity	&&\pageref{sec:mult}	\\
\ref{sec:counting}	&The abelian quotient theorem:			\\
			&\ \ proof by counting	&&\pageref{sec:counting}\\
\ref{sec:index}		&The abelian quotient theorem:			\\
			&\ \ proof by prime-index subgroups	
						&&\pageref{sec:index}	\\
\end{tabular}
\end{center}

\end{titlepage}


\section*{Introduction}	

Perfect numbers are an ancient object of study. A number is called perfect if
it is the sum of its proper divisors---for instance, 6 is perfect, since $6
= 1+2+3$. It is straightforward to classify the even perfect numbers, but it
is a long-standing question as to whether there are any odd perfect numbers
at all.

This article generalizes the notion of `perfection' from numbers to
groups. We define what it means for a group to be perfect, explain in what
sense this is a generalization of the notion for numbers, and go on to give
some theory and examples of perfect groups. Signposts are provided for the
reader not well versed in group theory, so that at least the rough shape of
the ideas should be discernible.

The first properties of perfect numbers are summarized in
Section~\ref{sec:numbers}. In Section~\ref{sec:basics} we give the definition
of a perfect group and look at some examples. Section~\ref{sec:mult} is
devoted to `multiplicativity'. This shows some close parallels with the world
of numbers, including the results of Section~\ref{sec:numbers}, and the new
theory also enables us to give some more interesting examples of perfect
groups than was possible previously. The climax of the article, such as it
is, is a theorem concerning the abelian quotients of perfect groups, a
corollary of which classifies the perfect abelian groups. Two rather
different proofs of this result are offered, one in each of
Sections~\ref{sec:counting} and~\ref{sec:index}.

\paragraph*{Acknowledgements and Apologies}

I would like to thank Robin Bhattacharyya for his careful reading-through of
an early version of this document, and Vin de Silva for reading a later
version.  Alan Bain made some useful suggestions on adapting it to an
undergraduate readership, and Simon Norton made some further helpful
remarks.  

I must also make two apologies.  First of all, I crave the reader's
indulgence for the use of the term `perfect group' when it is firmly
established to mean something else.  Faced with a group-theoretic concept
generalizing that of perfect number, any other name seemed unnatural.  My
second apology is for the lack of pointers to the literature: some of the
results included here are surely widely known, but I am not well enough
educated to provide references.

\sect{Perfect Numbers}	\label{sec:numbers}

Here we go over the basic properties of perfect numbers.

For any number $n$, define $D(n) = \sum_{d|n} d$, the sum of the divisors of
$n$, and call $n$ \defemph{perfect} if $D(n) = 2n$. By a `number' I mean a
positive integer.

\paragraph{Multiplicativity}	\label{para:nomult}
The function $D$ is \defemph{multiplicative}: that is to say, if $m_1$ and
$m_2$ are coprime (have no common divisors other than $1$) then $D(m_1
m_2) = D(m_1) D(m_2)$. To see this, first observe that any divisor $d$ of
$m_1 m_2$ can be written uniquely as $d_1 d_2$ where $d_i$ is a divisor of
$m_i$ ($i = 1, 2$); conversely, if $d_i$ is a divisor of $m_i$ ($i = 1, 2$)
then $d_1 d_2$ is a divisor of $m_1 m_2$.  Hence
\begin{eqnarray*}
D(m_1 m_2)	&=&\sum_{d|m_1 m_2} d			\\
		&=&\sum_{d_1|m_1,\, d_2|m_2} d_1 d_2	\\
		&=&\left(\sum_{d_1|m_1}d_1\right) 
		  \left(\sum_{d_2|m_2}d_2\right) 	\\
		&=&D(m_1)D(m_2),
\end{eqnarray*}
as required.

\paragraph{Even Perfect Numbers}
It is easy to classify the even perfect numbers: they are precisely those
numbers $2^{r-1}(2^{r}-1)$ where $r\geq 2$ and $2^{r}-1$ is prime. (Of
course, computing which values of $r$ make $2^{r}-1$ prime is itself a hard
problem.) The first three perfect numbers are $2\times 3 = 6$, $4\times 7 =
28$, and $16 \times 31 = 496$.

In one direction, suppose that $r\geq 2$ and $2^{r}-1$ is prime: then
\begin{eqnarray*}
D(2^{r-1}(2^r-1))&=&D(2^{r-1})D(2^r-1)
		   \makebox[0em]{\hspace*{21.5em}(by \ref{para:nomult})}\\
		&=&(1+2+2^2+ \cdots +2^{r-1})(1+2^r-1)	\\
		&=&(2^r-1)2^r				\\
		&=&2[2^{r-1}(2^r-1)],
\end{eqnarray*}
so $2^{r-1}(2^{r}-1)$ is an even perfect number.

In the other direction, suppose that $n$ is an even perfect number. Write
$n=2^{s}m$ where $s\geq 1$ and $m$ is odd: then $n$ being perfect says that
\[
D(2^{s}m) = 2\times 2^{s}m,
\]
i.e.
\[
(2^{s+1}-1)D(m) = 2^{s+1}m,
\]
i.e.
\begin{equation}
(2^{s+1}-1)(D(m)-m) = m.	\label{eq:evenperf}
\end{equation}
Hence $D(m)-m$ is a divisor of $m$, and since 
\[
2^{s+1}-1 > 2^{0+1}-1 =1,	
\]
it is a proper divisor of $m$. But $D(m)-m$ is by definition the sum of the
proper divisors of $m$, so $D(m)-m$ is the unique proper divisor of $m$. Thus
$m$ is prime and $D(m)-m = 1$, and by (\ref{eq:evenperf}), the latter means
that $m=2^{s+1}-1$. So $n=2^s(2^{s+1}-1)$ with $s\geq 1$ and $2^{s+1}-1$
prime, as required.

\sect{Definition and First Examples of Perfect Groups}
\label{sec:basics} 

In this section we define the notion of a perfect group, and search for
examples among some of the well-known families of groups (symmetric,
alternating, \ldots).  In fact, the only examples of perfect groups we will
find are cyclic, although by Section~\ref{sec:mult} we will have developed
enough theory to be able to exhibit some more interesting examples.

Of the examples below, only the cyclic groups (\ref{eg:cyclic}) and the
symmetric and alternating groups (\ref{eg:altsym}) will be needed later on.

The reader is reminded that a \defemph{normal subgroup} of a group $G$ is a
subset of $G$ which is the kernel of some homomorphism from $G$ to some other
group; equivalently, it is a subgroup $N$ of $G$ such that $gng^{-1} \in N$
for all $n\in N$ and $g\in G$.  We write $N\nsub G$ to mean that $N$ is a
normal subgroup of $G$.  From here on, `group' will mean `finite group'.

If $G$ is a group, define $D(G) = \sum_{N\nsub G} |N|$, the sum of the orders
of the normal subgroups of $G$, and say that $G$ is \defemph{perfect} if
$D(G) = 2|G|$.

\paragraph{Example: Cyclic Groups}	\label{eg:cyclic}
Let $C_n$ be the cyclic group of order $n$. Then $C_n$ has one normal
subgroup of order $d$ for each divisor $d$ of $n$, and no others, so $D(C_n)
= D(n)$ and $C_n$ is perfect just when $n$ is perfect. Thus perfect groups
provide a generalization of the concept of perfect numbers, and $C_6$,
$C_{28}$ and $C_{496}$ are all perfect groups.

\paragraph{Example: Symmetric and Alternating Groups}
\label{eg:altsym}
None of the symmetric groups $S_n$ or alternating groups $A_n$ is
perfect. If $n\geq 5$ then $A_n$ is simple and the only normal subgroups of
$S_n$ are 1, $A_n$ and $S_n$, so $D(A_n)$ and $D(S_n)$ are too small. For
$n\leq 4$, we have
\[\begin{array}{lcl}
D(A_1)=1,		&	&D(S_1)=1,		\\
D(A_2)=1,		&	&D(S_2)=1+2=3,		\\
D(A_3)=1+3=4,		&	&D(S_3)=1+3+6=10,	\\
D(A_4)=1+4+12=17,	&	&D(S_4)=1+4+12+24=41.
\end{array}\]

\paragraph{Example: $p$-Groups}
A (finite) \defemph{$p$-group} is a group of order $p^r$, where $p$ is prime
and $r\geq 0$. Lagrange's Theorem says that the order of any subgroup of a
group divides the order of the group, so if $G$ is a $p$-group then $D(G)
\equiv 1 \pmod{p}$. Hence no $p$-group is perfect.

\paragraph{Example: Dihedral Groups}	\label{eg:dihedral}
Let $E_{2n}$ be the dihedral group of order $2n$: that is, the group of all
isometries of a regular $n$-sided polygon. Of the $2n$ isometries, $n$ are
rotations (forming a cyclic subgroup of order $n$) and $n$ are
reflections. We examine the cases of $n$ odd and $n$ even separately.

\textbf{$n$ odd:} All reflections are in an axis passing through a
vertex and the midpoint of the opposite side, and any reflection is conjugate
to any other by a suitable rotation. Thus if $N\nsub E_{2n}$ and $N$ contains
a reflection, then $N$ contains all reflections; but $1\in N$ too, so $|N|
\geq n+1$, so $N=E_{2n}$. So any proper normal subgroup is inside the
rotation group $C_n$; conversely, any (normal) subgroup of $C_n$ is normal in
$E_{2n}$. Thus
\[
D(E_{2n}) = D(C_n) + 2n,
\]
and $E_{2n}$ is perfect if and only if $n$ is a perfect number.

\textbf{$n$ even:} The reflections split into two conjugacy classes, $R_1$
and $R_2$, each of size $n/2$: those in an axis through two opposite
vertices, and those in an axis through the midpoints of two opposite
sides. Write $C_{n/2}$ for the group of rotations by $2$ or $4$ or \ldots\ or
$n$ vertices, a subgroup of $E_{2n}$ which is cyclic of order $n/2$. Then we
can show that the smallest subgroup of $E_{2n}$ containing $R_i$ is $R_i \cup
C_{n/2}$, for $i=1$ and 2. Moreover, $R_i \cup C_{n/2}$ is of order $n$,
i.e.\ index 2, therefore normal in $E_{2n}$. So we have two different normal
subgroups, $R_1 \cup C_{n/2}$ and $R_2 \cup C_{n/2}$, of order $n$. We also
have the normal subgroups $\{1\}$ and $E_{2n}$, hence
\[
D(E_{2n}) \geq 1+n+n+2n > 4n
\]
and $E_{2n}$ is not perfect.

In summary, the perfect dihedral groups are in one-to-one correspondence with
the odd perfect numbers---so it is an open question as to whether there are
any.

\sect{Multiplicativity}	\label{sec:mult}

We proved in~\ref{para:nomult} that the function $D(n)$, on numbers $n$, was
multiplicative. The aim of this section is to prove an analogous result for
groups, and then to give some examples of nonabelian perfect groups by using
this result.

Some difficulties are present for the reader not acquainted with composition
series and the Jordan-H\"older Theorem. However, it is still possible for him
or her to understand an example (\ref{eg:s3xc5}) of a nonabelian perfect
group, provided that the following fact is taken on trust: if $G_1$ and $G_2$
are groups whose orders are coprime, and $G_1 \times G_2$ their direct
product, then $D(G_1 \times G_2) = D(G_1) D(G_2)$.  This done, the reader
may proceed to~\ref{eg:s3xc5} straight away.

The Jordan-H\"older Theorem states that any two composition series for a
group $G$ have the same set-with-multiplicities of factors, up to isomorphism
of the factors. I shall write this set-with-multiplicities as $c(G)$, and use
$+$ to denote the disjoint union (or `union counting multiplicities') of two
sets-with-multiplicities. Thus if 
\[
c(G) = \{C_2, C_2, C_5\} \mbox{\ \ and\ \ } c(H) = \{C_2, A_6\}
\]
then
\[
c(G)+c(H) = \{C_2, C_2, C_2, C_5, A_6\}.
\]
We will use the fundamental fact that if $K\nsub X$ then
\[
c(X) = c(X/K) + c(K).
\]

A pair of groups will be called \defemph{coprime} if they have no composition
factor in common; alternatively, we will say that one group is \defemph{prime
to} the other.  (In particular, if two groups have coprime orders then they
are coprime.) We will prove that $D$ is \defemph{multiplicative}: that is, if
$G_1$ and $G_2$ are coprime then $D(G_1 \times G_2) = D(G_1) D(G_2)$. First
of all we establish the group-theoretic analogue of a number-theoretic result
from Section~\ref{sec:numbers}---namely, the second sentence
of~\ref{para:nomult}. 

\paragraph{Proposition}	\label{para:nsubprod}
\textit{%
Let $G_1$ and $G_2$ be coprime groups. Then the normal subgroups of $G_1
\times G_2$ are exactly the subgroups of the form $N_1 \times N_2$, with $N_1
\nsub G_1$ and $N_2 \nsub G_2$.}

\pf If $N_1 \nsub G_1$ and $N_2 \nsub G_2$ then $N_1 \times N_2 \nsub G_1
\times G_2$; conversely, suppose $N \nsub G_1 \times G_2$. Write $\pi_i: G_1
\times G_2 \go G_i$ ($i=1$, 2) for the projections, and regard $G_1$ as a
normal subgroup of $G_1 \times G_2$ by identifying it with $G_1 \times
\{1\}$, and similarly $G_2$.  We have
\[
\pi_1 N \iso \frac{N}{\ker (\pi_1 |_N)} = \frac{N}{G_2 \cap N},
\]
so by the `fundamental fact' above,
\[
c(N) = c(\pi_1 N) + c(G_2 \cap N);
\]
and therefore by symmetry
\[
c(\pi_1 N) + c(G_2 \cap N) = c(\pi_2 N) + c(G_1 \cap N).
\]
But $c(\pi_i N) \sub c(G_i)$ and $G_1$ and $G_2$ are coprime, so $c(\pi_1 N)$
and $c(\pi_2 N)$ have no element in common; similarly $c(G_i \cap N) \sub
c(G_i)$, so $c(G_1 \cap N)$ and $c(G_2 \cap N)$ have no element in
common. Hence $c(\pi_i N) = c(G_i \cap N)$. We also know that $c(X)$
determines the order of a group $X$ and that $G_i \cap N \sub \pi_i N$, so
in fact $G_i \cap N = \pi_i N$. Thus
\[
\pi_1 N \times \pi_2 N = (G_1 \cap N) \times (G_2 \cap N) \sub N,
\]
and as always
\[
N \sub \pi_1 N \times \pi_2 N,
\]
so $N = \pi_1 N \times \pi_2 N$, with $\pi_i N \nsub G_i$.
\done

\paragraph{Corollary}
\textit{$D$ is multiplicative.}

\pf
\marginpar{\hspace{-1.8em}\raisebox{-22.4ex}{$\Box$}} 
This is a direct analogue of~\ref{para:nomult}. For by~\ref{para:nsubprod},
if $G_1$ and $G_2$ are coprime then
\begin{eqnarray*}
D(G_1 \times G_2)	&=&\sum_{N_1 \nsub G_1,\, N_2 \nsub G_2}
				|N_1 \times N_2|	\\
			&=&\sum_{N_1 \nsub G_1} \sum_{N_2 \nsub G_2}
				|N_1| |N_2|	\\
			&=&D(G_1) D(G_2).
\end{eqnarray*}

\paragraph*{}
We can now exhibit three nonabelian perfect groups.

\paragraph{Example: $S_3 \times C_5$}	\label{eg:s3xc5}
The group $S_3 \times C_5$, of order 30, is perfect. For $S_3$ and $C_5$ have
coprime orders (6 and 5), so are coprime, so 
\begin{eqnarray*}
D(S_3 \times C_5)	&=&D(S_3) D(C_5)	\\
			&=&(1+3+6) \times (1+5)	\\
			&=&60			\\
			&=&2|S_3 \times C_5|.
\end{eqnarray*}

\paragraph{Example: $A_5 \times C_{15128}$}
We present this example (of order $907\,680$) along with the method by which
it was found. Firstly, $A_5$ is a simple group of order $5!/2 = 60$. Now, let
us try to find a perfect group $G$ of the form $G = A_5 \times G_1$ where
$G_1$ is some group prime to $A_5$. Since
\[
D(A_5)/|A_5| = 61/60,
\]
we need to find a $G_1$ such that
\[
D(G_1)/|G_1| = 120/61.
\]
Let us look for such a group $G_1$ amongst those of the form $G_1 = C_{61}
\times G_2$, where $G_2$ is prime to $C_{61}$ and $A_5$. Since
\[
D(C_{61})/|C_{61}| = 62/61,
\]
we need to find a $G_2$ such that
\[
D(G_2)/|G_2| = 120/62=60/31.
\]
In turn, let us look for such a group $G_2$ amongst those of the form $G_2 =
C_{31} \times G_3$, where $G_3$ is prime to $C_{31}$, $C_{61}$ and
$A_5$. Since
\[
D(C_{31})/|C_{31}| = 32/31,
\]
we need to find a $G_3$ such that
\[
D(G_3)/|G_3| = 60/32 = 15/8.
\]
This is satisfied by $G_3 = C_8$, and the groups $A_5$, $C_{61}$,
$C_{31}$ and $C_8$ are pairwise coprime.  Thus if
\begin{eqnarray*}
G	&=&A_5 \times C_{61} \times C_{31} \times C_8	\\
	&=&A_5 \times C_{61\times 31\times 8}		\\
	&=&A_5 \times C_{15128}
\end{eqnarray*}
then $G$ is perfect.

\paragraph{Example: $A_6 \times C_{366776}$} 
By the same technique we get this next example, of order $132\,039\,360$.
This time, we start with the simple group $A_6$ of order $6!/2 = 360$, and
the sequence of groups $A_6$, $C_{361}$, $C_{127}$, $C_8$ `works' in the
sense of the previous example. The details are left to the reader; note that
$361 = 19^2$ and that 127 is prime.

\sect{\sloppy The Abelian Quotient Theorem: Proof by Counting}
\label{sec:counting}

In each of the next two sections we present a separate proof of our main
classification result, the abelian quotient theorem. The two proofs have
rather different flavours, and each produces its own insights, which is why
both are included. We start with the more elementary of the two.

An \defemph{abelian quotient} of a group $G$ is just a quotient of $G$ which
is abelian. That is, it's an abelian group $A$ for which there exists a
surjective homomorphism $G \go A$; alternatively, it's an abelian group
isomorphic to $G/K$ for some normal subgroup $K$ of $G$.  We will prove:

\paragraph{Abelian Quotient Theorem}	\label{thm:abqt}
\begin{quote}
\textit{%
If $G$ is a group with $D(G)\leq 2|G|$ then any abelian quotient of $G$ is
cyclic.}
\end{quote}

\paragraph*{}This result has the following corollaries, the second of which
says that abelian perfect groups `are' just perfect numbers: 

\paragraph{Corollaries}	\label{cors:abqt}
{\itshape
\begin{enumerate}
\item \label{cor:abqtpfcy} 
If $G$ is a perfect group then any abelian quotient of $G$ is cyclic.
\item \label{cor:abpfcy}
The perfect abelian groups are precisely the cyclic groups $C_n$ of order
$n$ with $n$ perfect.
\end{enumerate}}

\pf 
Part (\ref{cor:abqtpfcy}) is immediate. For (\ref{cor:abpfcy}), if $A$ is
perfect abelian then $A$ is an abelian quotient of the perfect group $A$,
hence $A$ is cyclic. But we have already seen~(\ref{eg:cyclic}) that the
perfect cyclic groups correspond exactly to the perfect numbers.  \done

\paragraph*{}

(Those who know about such things will recognize that the theorem could be
stated more compactly in this way: if $G$ is a group with $D(G) \leq 2|G|$
then $G^\mathrm{ab}$ is cyclic. Here $G^\mathrm{ab}$ is the
\defemph{abelianization} of $G$: it is an abelian quotient of $G$ with the
property that any abelian quotient of $G$ is also a quotient of
$G^\mathrm{ab}$. In particular, if $A$ is abelian then $A^\mathrm{ab} \iso
A$, which is how we would deduce Corollary~\ref{cors:abqt}(\ref{cor:abpfcy})
from this formulation.)

The proof of the abelian quotient theorem given in this section uses two
ingredients. The first is a new way of evaluating $D(G)$:

\paragraph{Lemma}	\label{lemma:count}
\textit{%
For any group $G$,
\[
D(G) = \sum_{g\in G} |\{\mbox{\rm normal subgroups of $G$ containing $g$}\}|.
\]}

\pf
We have
\begin{eqnarray*}
D(G)	&=&\sum_{N\nsub G} |N|			\\
	&=&|\{(N,g): N\nsub G,\,g\in N\}|	\\
	&=&\sum_{g\in G} |\{N: N\nsub G,\,g\in N\}|.
\end{eqnarray*} \done

\paragraph*{}
The second ingredient is the `standard' fact that the inverse image (under a
homomorphism) of a normal subgroup is a normal subgroup. For let $\pi: G_1
\go G_2$ be a homomorphism of groups, and let $N\nsub G_2$. Then $N$ is the
kernel of the natural homomorphism $\phi: G_2 \go G_2 / N$, in other words,
$N=\phi^{-1}\{0\}$.  So 
\[
\pi^{-1}N = 
\pi^{-1}\phi^{-1}\{0\} =
(\phi\of\pi)^{-1}\{0\},
\]
i.e.\ $\pi^{-1}N$ is the kernel of the homomorphism $\phi\of\pi: G_1 \go G_2
/ N$.  Thus $\pi^{-1}N$ is a normal subgroup of $G_1$.

We are now ready to assemble these ingredients into the following
proposition, from which the abelian quotient theorem follows immediately.
Two pieces of terminology will be used.  An element $h$ of $G$ is called a
\defemph{normal generator} of $G$ if the only normal subgroup of $G$
containing $h$ is $G$ itself.  A group is called \defemph{simple} if it has
precisely two normal subgroups---inevitably, the whole group and the
one-element subgroup.

\paragraph{Proposition}
{\itshape
Let $G$ be a group. 
\begin{enumerate}
\item 	\label{prop:deficient}
If $D(G) \leq 2|G|$ then $G$ has a normal generator.
\item If $G$ has a normal generator then any abelian quotient of $G$ is
cyclic.
\end{enumerate}}

\pf
\begin{enumerate}
\item
By Lemma~\ref{lemma:count}, $D(G) \leq 2|G|$ if and only if the mean over all
$g\in G$ of
\[
\nu(g) := |\{\mbox{\rm normal subgroups of $G$ containing $g$}\}|
\]
is $\leq 2$. If $G$ is not simple or trivial then $\nu(1_G)\geq 3$ (where
$1_G$ is the identity element of $G$); so for the mean to be $\leq 2$,
there must be some $h\in G$ for which $\nu(h)=1$---and this says exactly
that $h$ is a normal generator of $G$. On the other hand, if $G$ is simple
then any nonidentity element of $G$ is a normal generator, and if $G$ is
trivial then $1_G$ is a normal generator. So (\ref{prop:deficient}) is proved
in all cases.

\item
Let $A$ be an abelian quotient of $G$, with $\pi: G \go A$ a surjective
homomorphism, and let $h$ be a normal generator of $G$.  Then $\pi(h)$ is a
normal generator of $A$: for if $K \nsub A$ and $\pi(h) \in K$ then
$\pi^{-1}K$ is a normal subgroup of $G$ containing $h$, so $\pi^{-1}K = G$;
and since $\pi$ is surjective, this means that $K=A$.  But $A$ is abelian,
so all subgroups are normal, so the fact that $\pi(h)$ is a normal generator
of $A$ says that the only subgroup of $A$ containing $\pi(h)$ is $A$ itself.
And this in turn says exactly that the cyclic subgroup generated by $\pi(h)$
is $A$ itself.  
\done
\end{enumerate} 

\sect{\sloppy The Abelian Quotient Theorem: Proof by Prime-Index Subgroups}
\label{sec:index}

This last section is devoted to a second proof of the abelian quotient
theorem,~\ref{thm:abqt}.  This time, the proof reveals something about the
normal subgroup structure of a perfect group $G$: namely, that $G$ has at
most one normal subgroup of each prime
index~(\ref{propn:index}(\ref{prop:pfone})). It is a corollary of this that
any abelian quotient of $G$ is cyclic.

This section assumes some more sophisticated group theory than the last.

\paragraph{Lemma}	\label{lemma:index}
\textit{%
Let $G$ be a group and $p$ a prime: then the number of normal subgroups of
$G$ with index $p$ is
\[
\frac{p^r-1}{p-1} = 1+p+\cdots +p^{r-1},
\]
for some $r\geq 0$.}

\paragraph*{Remark} `Usually' $r=0$, in which case both sides of the equation
evaluate to~$0$.

\pf

For this proof we write the cyclic group of order $p$ additively, as
$\zed/p\zed$.  We also write $\mathrm{Hom}(G,\zed/p\zed)$ for the set of all
homomorphisms $G \go \zed/p\zed$, and $\mathrm{Aut}(\zed/p\zed)$ for the set
of all automorphisms of the group $\zed/p\zed$ (that is, invertible
homomorphisms $\zed/p\zed \go \zed/p\zed$).

The key observation is that a normal subgroup of $G$ of index $p$ is just the
kernel of a surjection from $G$ to $\zed/p\zed$.

All but one element of $\mathrm{Hom}(G,\zed/p\zed)$ is surjective, and the
remaining one is trivial. Two surjections $\pi, \phi: G \go \zed/p\zed$ have
the same kernel if and only if $\pi = \alpha\of\phi$ for some $\alpha\in
\mathrm{Aut}(\zed/p\zed)$; moreover, if such an $\alpha$ exists for $\pi$ and
$\phi$ then it is unique.  So the nontrivial elements of
$\mathrm{Hom}(G,\zed/p\zed)$ have
\[
\frac{|\mathrm{Hom}(G,\zed/p\zed)| - 1}{|\mathrm{Aut}(\zed/p\zed)|}
\]
different kernels between them. In other words, there are this many index-$p$
normal subgroups of $G$. We now just have to evaluate
$|\mathrm{Hom}(G,\zed/p\zed)|$ and $|\mathrm{Aut}(\zed/p\zed)|$.

Firstly, $\zed/p\zed$ is cyclic with $p-1$ generators, so
$|\mathrm{Aut}(\zed/p\zed)| = p-1$.

Secondly, $\zed/p\zed$ is abelian, so $\mathrm{Hom}(G,\zed/p\zed)$ forms an
abelian group under pointwise addition. Each element has order 1 or $p$, so
$\mathrm{Hom}(G,\zed/p\zed)$ can be given scalar multiplication over the
field $\zed/p\zed$, and thus becomes a finite vector space over
$\zed/p\zed$. This vector space has a dimension $r\geq 0$, and then
$|\mathrm{Hom}(G,\zed/p\zed)| = p^r$. (Alternatively, Cauchy's Theorem gives
this result.) 

The lemma is now proved.  
\done

\paragraph*{}

Let us temporarily call a group $G$ \defemph{tight} if for each prime $p$,
$G$ has at most one normal subgroup of index $p$.  Putting together the three
parts of the following proposition gives us our second proof of the abelian
quotient theorem.

\paragraph{Proposition}		\label{propn:index}
{\itshape
\begin{enumerate}
\item	\label{prop:pfone}
A group $G$ with $D(G)\leq 2|G|$ is tight.
\item 	\label{prop:quotient}
A quotient of a tight group is tight. 
\item
A tight abelian group is cyclic.
\end{enumerate}}

\pf
\begin{enumerate}
\item
For each prime $p$, we have
\[
2|G| \geq D(G) \geq |G| + \frac{p^r-1}{p-1} \cdot \frac{|G|}{p},
\]
where $r$ is as in Lemma~\ref{lemma:index}. If $r\geq 2$ then 
\[
\frac{p^r-1}{p-1} \cdot \frac{|G|}{p} \geq (p+1)\cdot \frac{|G|}{p} > |G|,
\]
giving a contradiction. Thus $r$ is 0 or 1, and so $\frac{p^r-1}{p-1}$ is 0
or 1.

\item
Let $\pi: G_1 \go G_2$ be a surjective homomorphism. If $N$ and $N'$ are
distinct normal subgroups of $G_2$ with index $p$, then $\pi^{-1}N$ and
$\pi^{-1}N'$ are distinct normal subgroups of $G_1$ with index $p$.

\item
For this we invoke the classification theorem for finite abelian groups,
which tells us that for any abelian group $A$ there exist primes $p_1,
\ldots, p_n$ and numbers $t_1, \ldots, t_n \geq 1$ such that
\[
A \iso C_{p_1^{t_1}} \times \cdots \times C_{p_n^{t_n}}.
\]
Suppose that $p_i = p_j$ ($=p$, say) for some $i\neq j$. Then, since $t_i
\geq 1$, $C_{p^{t_i}}$ has a (normal) subgroup $N_i$ of index $p$; and
similarly $C_{p^{t_j}}$.  Hence $N_i \times C_{p^{t_j}}$ and $C_{p^{t_i}}
\times N_j$ are distinct index-$p$ subgroups of $C_{p^{t_i}} \times
C_{p^{t_j}}$, and $C_{p^{t_i}} \times C_{p^{t_j}}$ is not tight.  Since
$C_{p^{t_i}} \times C_{p^{t_j}}$ is a quotient of $A$,
part~(\ref{prop:quotient}) implies that $A$ is not tight either.  Thus if $A$
\emph{is} tight then all the $p_k$'s are distinct, so that
\[
A \iso C_{p_1^{t_1}p_2^{t_2} \cdots p_n^{t_n}}.
\]\done
\end{enumerate}

\paragraph*{}
There are still other lines of proof for the abelian quotient theorem. In
part~(\ref{prop:quotient}) of the Proposition, the fact that $p$ was prime
was quite irrelevant, and in just the same manner we can prove that
\[
\frac{D(G_1)}{|G_1|} \geq \frac{D(G_2)}{|G_2|}
\]
whenever $G_2$ is a quotient of $G_1$.  (If $\pi: G_1 \go G_2$ is the
quotient map, with kernel of order $k$, then a normal subgroup $N$ of $G_2$
gives rise to a normal subgroup $\pi^{-1}N$ of $G_1$ of order $k|N|$.)  Thus
if $G$ is a group with $D(G) \leq 2|G|$ and $A$ is an abelian quotient of $G$
then $D(A) \leq 2|A|$. So we have reduced the abelian quotient theorem to the
abelian case: if $A$ is abelian and $D(A) \leq 2|A|$ then $A$ is cyclic.
Certainly this is provable by methods derived from one of the two proofs of
the general case, but other approaches exist; I leave that for the reader.

\section*{Further Thoughts}

We finish with some general speculative thoughts, roughly in order of the
material above.

The chosen definition of the function $D$, and therefore of perfect group, is
one amongst many candidates.  We defined $D$ to be the sum of the orders of
the normal subgroups, but we could change `normal subgroups' to `subgroups',
`characteristic subgroups', `subnormal subgroups', \ldots, or we could define
$D$ to be the sum of the \emph{indices} of the normal subgroups, etc.  In all
cases we preserve the identity $D(C_n) = D(n)$, but only in some of them does
$D$ remain multiplicative (a feature we probably like).

More abstractly, this article was about lifting the classical function $D: \{
\mathrm{numbers} \} \go \{ \mathrm{numbers} \}$ to a function $D: \{
\mathrm{groups} \} \go \{ \mathrm{numbers} \}$.  We might consider it natural
to go the whole hog and create a function assigning not just a number, but
some kind of algebraic structure, to each group $G$.  I do not know of any
very useful way to do this.

In number theory there is a whole body of work on multiplicative functions of
integers, which include the number-of-divisors function, the sum-of-divisors
function, the Euler function $\phi$, and the M\"obius function $\mu$. In the
world of groups we have at least the beginning of an analogue. For let $F$ be
a multiplicative function from groups to numbers: then just as in
Corollary~3.2, the function $F': G \longmapsto \sum_{N\nsub G} F(N)$ is
multiplicative. For instance, if $F$ is the function with constant value~$1$
then $F'$ gives the number of normal subgroups of a group, and is
multiplicative.

The abelian quotient theorem says that if $D(G) \leq 2|G|$ then $G$ has some
special property expressible in standard group-theoretic terms. We can prove
this in at least two ways, but it seems rather more challenging to prove
something in the other direction: that if $D(G)$ is `too big' then $G$ has a
certain form.

Finally, we can make various conjectures on perfect groups, based on the
skimpy evidence above: for instance, `there are no odd-order perfect
groups', or `there are infinitely many nonabelian perfect groups'.
Example~\ref{eg:dihedral}, on the dihedral groups, tells us that classifying
the even-order perfect groups is at least as hard as determining whether
there are any odd perfect numbers.  Clearly such problems are unlikely to
be easy to solve.

\end{document}